\documentclass[preprint,11pt]{elsarticle}

\usepackage{amssymb}
\usepackage{mathrsfs}
\usepackage{amsthm}
\usepackage{amsmath}
\usepackage{epstopdf}
\usepackage{epsfig}
\usepackage{xcolor}
\usepackage{graphicx}
\usepackage{subfigure}
\usepackage{longtable,booktabs}
\usepackage{lineno,hyperref}
\usepackage{threeparttable}
\usepackage{mathtools}
 \usepackage{verbatim}

\newtheorem{theorem}{Theorem}[section]
\newtheorem{lemma}{Lemma}[section]

\newtheorem{proposition}{Proposition}[section]
\newtheorem{definition}{Definition}[section]

\newtheorem{example}{Example}[section]

\begin{document}

\begin{frontmatter}



\title{Classification of permanence and impermanence for a Lotka-Volterra model of three competing
species with seasonal succession}

\author[add1]{Lei Niu}\ead{lei.niu@dhu.edu.cn}


\author[add3]{Xizhuang Xie\corref{cor1}}\ead{xzx@hqu.edu.cn}

\address[add1]{Department of Mathematics, Donghua University, Shanghai
201620, China}
\address[add3]{School of Mathematical Sciences, Huaqiao University
Quanzhou, Fujian 362021, China}

\cortext[cor1]{Corresponding author}
\begin{abstract}
In this paper, we are concerned with the permanence of a Lotka-Volterra model of three competing species with seasonal
succession. Based on the existence of a carrying simplex, that is a globally attracting hypersurface of codimension one, we provide a complete classification of the permanence and impermanence in terms of inequalities on the parameters of this model. Moreover, invariant closed curves are found numerically in the permanent classes, which means that the positive fixed point of the Poincar\'e map in the permanent classes is not always globally asymptotically stable.
\end{abstract}


\end{frontmatter}



\section{Introduction}
A fundamental issue in population biology is to find the coexistence conditions for interacting species. Permanence is one of such concepts, which 
is designed to provide an answer to such questions as which species, in a mathematical model of interacting species, will survive over the 
long term. By establishing an initial-condition-independent positive lower bound for the long-term value of a component of a dynamical system
such as population size, the permanence theory can give a mathematically rigorous answer to the question of permanence. There has been much attention on Kolmogorov type differential equation models (see \cite{Hofbauer1981,Hutson1984,Hutson1988,Hofbauer1988}) and difference equation models (see \cite{Hutson1982,hofbauer1987coexistence,Lu1999,Kon2001,Kon2004}). It has also been characterized from a topological perspective (see \cite{Butler1986,Hofbauer1989,Garay1989,Garay2003,Gyllenberg2020b}).
However, the important case of models of differential equations with seasonal succession (see, for example, \cite{Hsu-Zhao,k2010,k2012,k2013,lk2011,tang2017,xiao2016,xie2021,zhang-zhao2013}) appears to have received relatively little attention from this point of view. As we know, interactive species often live in
a fluctuating environment. Due to the seasonal or daily variations, populations experience a periodic dynamical environment driven by both internal dynamics of species interactions and external forcing (\cite{yohe2003,sommer1986,walther2002}). The vector fields of the models with seasonal succession are usually discontinuous and periodic with respect to time. For such models, mathematical obstacles include the existence and number of the positive periodic solutions, and estimates
for the Floquet multipliers, etc..  To our best knowledge, even for two
species models, progress has been somewhat limited.

 In this paper, we shall focus on the Lotka-Volterra  model of  three competing species  with seasonal
succession
\begin{equation}\label{seasonal-system}
 \begin{dcases}
\frac{d x_{i}}{d t}=-\mu_{i} x_{i}, \quad t\in[k \omega,  k \omega+(1-\varphi) \omega), \\
\frac{d x_{i}}{d t}=x_{i}\left(b_{i}-\sum_{j=1}^{3} a_{i j} x_{j}\right),\quad t\in [k \omega+(1-\varphi) \omega,  (k+1) \omega),
\end{dcases}
\end{equation}
where $k \in \mathbb{Z}_{+}, \varphi \in (0,1]$ and $\omega$, $\mu_{i}$, $b_{i}$, and $a_{i j}$, $i,j=1,2,3$, are all positive constants. The model \eqref{seasonal-system} is a time-periodic system in a seasonal succession environment, where the overall period is $\omega$, and $\varphi$ stands for the switching proportion of a period between the linear system
\begin{equation}\label{linear-system}
    \frac{d x_{i}}{d t}=-\mu_{i} x_{i}, ~i=1, 2, 3
\end{equation}
and the classical Lotka-Volterra competition model
\begin{equation}\label{LV-system}
    \frac{d x_i}{dt}=x_i\left(b_i- \sum_{j=1}^3 a_{ij}x_j \right),\quad i=1,2,3.
\end{equation}
Biologically, succession can be thought of as a series of transitions between systems (\ref{linear-system}) and (\ref{LV-system}) where
$\varphi$ is used to describe the proportion of the period in the good season in which the species follow the Lotka-Volterra system \eqref{LV-system}, while $(1-\varphi)$ represents the proportion of the period in the bad season in which the species die exponentially according to system \eqref{linear-system}.

For the Lotka-Volterra  model of two competing species  with seasonal
succession, it was thoroughly analyzed by Hsu and Zhao \cite{Hsu-Zhao} and Niu et al. \cite{Niu-Wang-Xie}. Hsu and Zhao \cite{Hsu-Zhao} showed that every solution of the two species model is asymptotic to some periodic solution and there are a total of four dynamical scenarios. Niu et al. \cite{Niu-Wang-Xie} proved that the associated Poincar\'e map of the model admits a globally attracting hypersurface of codimension one, called the carrying simplex, which carries all the relevant long-term dynamics. We refer the readers to, for example,  \cite{diekmann2008carrying,Gyllenberg2020b,hirsch1988,hirsch2008existence,LG,jiang2015,smith1986,ruiz2011exclusion,wang2002uniqueness,Z993} for a detailed discussion on the theory and technique of the carrying simplex. Based on the existence of a carrying simplex, the authors recently classified the dynamics of the Poincar\'{e} map associated with the three species model \eqref{seasonal-system} into 33 equivalence classes by which the existence and multiplicity of the positive periodic solution for each class was presented.

Our object in this paper is to establish the criteria of permanence and impermanence for the model
\eqref{seasonal-system}. The main mathematical tools involve the theories of discrete dynamical systems and the carrying simplex. We will consider the discrete dynamical system defined by the iterates of the Poincar\'{e} map associated with the model \eqref{seasonal-system}. We give a complete classification of permanence and impermanence in terms of inequalities on parameters of this model based on the 33 equivalence classes provided in \cite{NWX2023}. Our numerical simulations show that invariant closed curves can occur in some permanent classes, which means that the positive fixed point of the Poincar\'e map in the permanent classes is not always globally asymptotically stable.


\section{Preliminaries}\label{section:3}
We use $\mathbb{R}^3_+$ to denote the nonnegative cone $\{x\in \mathbb{R}^3: x_i\geq 0,~ i=1,2,3\}$. The interior of $\mathbb{R}^3_+$ is the open cone $\dot{\mathbb{R}}^3_+:= \{x\in \mathbb{R}^3_+: x_i>0,~i=1,2,3\}$ and the boundary of $\mathbb{R}^3_+$ is $\partial \mathbb{R}^3_+:=\mathbb{R}^3_+\setminus \dot{\mathbb{R}}^3_+$. We denote by $\mathbb{H}^+_{i}$ the $i$-th positive coordinate axis and by ${\Pi}_i=\{x\in \mathbb{R}^3_+:x_i=0\}$ the $i$-th coordinate plane, $i=1,2,3$. We say $x$ is nonnegative if it belongs to $\mathbb{R}^3_+$ while positive if it belongs to $\dot{\mathbb{R}}^3_+$.  The symbol $0$  stands for both the origin of $\mathbb{R}^3$ and the real number $0$.

The Cauchy problem for \eqref{seasonal-system} has a unique global solution whenever the initial data $x(t_0)=x_0$ belongs to $\mathbb{R}_+^3$. We then may introduce, for any $t \geq t_0$, the solution operator $\Psi(t, t_0)$: it is defined on an open set $W\subset \mathbb{R}^3$ containing $\mathbb{R}_+^3$, and maps the initial datum $x_0$ into the solution $x(t)$ at time $t$.

It is easy to see that $\Psi$ has the properties:
\begin{itemize}
    \item[(i)] $\Psi$ is continuously
differentiable with respect to the initial data;
\item[(ii)] $\Psi(t, t_0)\mathbb{H}^+_{i} \subset \mathbb{H}^+_{i}$, $\Psi(t, t_0){\Pi}_i \subset {\Pi}_i$, $\Psi(t, t_0)\dot{\mathbb{R}}_{+}^3 \subset \dot{\mathbb{R}}_{+}^3$,\quad $\forall t \geq t_0$, $i=1,2,3$.
\item[(iii)] $ \Psi(t, s)\circ \Psi(s, t_0)=\Psi(t, t_0)$, \quad $\forall t \geq s \geq t_0$,  $t_0 \in \mathbb{R}$.
\item[(iv)] $\Psi(t+\omega, t_0+\omega)=\Psi(t, t_0)$, \quad $\forall t \geq t_0$,  $t_0 \in \mathbb{R}$.
\end{itemize}
When $t_0=0$, we write $\Psi(t, x)$, for $t\geq 0$, to denote the solution of \eqref{seasonal-system} with the initial data $\Psi(0)=x$. We are interested in the dynamics of $\Psi(t, x)$ in $\mathbb{R}_+^3$.

Consider the associated Poincar\'{e} map $\mathcal{P}$ of $\Psi(t, x)$ with
$\mathcal{P}(x)=\Psi(\omega, x)$, $x\in W$. Let $L:\mathbb{R}^3\to \mathbb{R}^3$ be the linear map
\begin{equation}\label{linear-map}
    x\mapsto (e^{-\mu_{1}(1-\varphi) \omega} x_{1}, e^{-\mu_{2}(1-\varphi) \omega} x_{2}, e^{-\mu_{3}(1-\varphi) \omega} x_{3})
\end{equation}
with $x \in \mathbb{R}^{3}$.
We denote by $\Phi_{t}(x)$ the solution map associated with the Lotka-Volterra competitive system \eqref{LV-system}. Then, we have
$$
\mathcal{P}(x)=\Phi_{\varphi \omega}\left(L x\right), \quad x \in W.
$$
Obviously, $\mathcal{P}(\mathbb{R}_+^3)\subseteq \mathbb{R}_+^3$. Since the coordinate axes and planes are forward invariant under $\mathcal{P}$ by the properties of $\Psi$, we can rewrite the map $\mathcal{P}$ as:
\begin{equation}\label{equ:Kol}
    \mathcal{P}(x_{1}, x_2, x_3)=\left(x_{1} F_{1}(x), x_{2} F_{2}(x), x_{3} F_{3}(x)\right),
\end{equation}
where
$$
F_{i}(x):= \begin{cases}\frac{\mathcal{P}_{i}(x)}{x_{i}} & \text { if } x_{i} \neq 0, \\[3pt]
\frac{\partial \mathcal{P}_{i}}{\partial x_{i}}(x) & \text { if } x_{i}=0.\end{cases}
$$
It is clear that $F_i(x)$ are continuous with respect to $x\in \mathbb{R}_+^3$, $i=1,2,3$.
\begin{lemma}\label{lemma:Fi}
     $F_i(x)>0$ for all $x\in \mathbb{R}_{+}^{3}$, $i=1,2,3$.
\end{lemma}
\begin{proof}
   First, note that by the property (ii) of $\Psi$, one has $\mathcal{P}_i(x)>0$ if $x_i>0$, and hence $F_i(x)>0$ if $x_i>0$. Now it suffices to prove that $F_i(x)>0$ if $x_i=0$. We may assume, without loss of generality, that $x_3=0$. By \eqref{equ:Kol}, one has
$$
\frac{\partial \mathcal{P}_{3}}{\partial x_{3}}(x)=F_{3}(x)\  \
{\rm and} \   \  \
\frac{\partial \mathcal{P}_{3}}{\partial x_{j}}(x)=0, \quad \forall j \neq 3.
$$
Let $g_{i}(x)=b_{i}-\sum_{j=1}^{3} a_{i j} x_{j}$ and $f(x)=(f_1(x),f_2(x),f_3(x))$ with $f_i(x)=x_ig_i(x)$, $i=1,2, 3$.
Then one has
\begin{equation}
D f(x)=\left[\begin{array}{cccc}
g_{1}(x)-a_{11} x_{1} & -a_{12} x_{1} & -a_{1 3} x_{1} \\
-a_{21} x_{2} & g_{2}(x)-a_{22} x_{2} &-a_{2 3} x_{2} \\
-a_{3 1} x_{3} & -a_{3 2} x_{3} & g_{3}(x)-a_{33}x_3
\end{array}\right].
\end{equation}
Let $W(t,x)=D_{x} \Phi_{t}(x)$ and $U(t,x)=D f(\Phi_{t}(x))$. Then
\begin{equation}\label{equ-DPhi}
    \frac{d W(t)}{d t}=U(t,x) \cdot W(t), \quad W(0)=I.
\end{equation}
Since $\mathcal{P}(x)=\Phi_{\varphi \omega} (Lx)$, one has
\begin{equation}\label{equ-DPx}
    D \mathcal{P}(x)=W(\varphi \omega,L x) \cdot D L,
\end{equation}
where $DL=\mathrm{diag}\left(e^{-\mu_{1}(1-\varphi) \omega}, e^{-\mu_{2}(1-\varphi) \omega}, e^{-\mu_{3}(1-\varphi) \omega}\right)$. Let $v(t,x):=\Phi_{t} (x)$. Then $v_3(t,Lx)=0$ by the expression \eqref{LV-system} and $x_3=0$. Thus, $$U(t,Lx)_{(3,3)}=g_3(\Phi_{t} (Lx))\quad \text{and}\quad U(t,Lx)_{(3,j)}=0~ \mathrm{for~} j\neq 3.$$
It follows that
\begin{equation}
    W(\varphi\omega,Lx)_{(3,3)}=\exp\left\{\int_0^{\varphi\omega} g_3(\Phi_{t} (Lx))dt\right\}, 
\end{equation}
$$W(\varphi\omega,Lx)_{(3,j)}=0~ \mathrm{for~} j\neq 3.$$
Therefore, by \eqref{equ-DPx}, one has
$$F_{3}(x)=\exp\left\{\int_0^{\varphi\omega} g_3(\Phi_{t} (Lx))dt-\mu_{3}(1-\varphi) \omega\right\}>0.$$
\end{proof}

\begin{definition}[\cite{Hofbauer1988, Kon2004, Smith2010}]
The system \eqref{seasonal-system} is said to be permanent if there are numbers $\delta,M>0$ such that
$$
\delta\leq \liminf_{t\to +\infty}\Psi_i(t, x)\leq \limsup_{t\to +\infty} \Psi_i(t, x) \leq M,\quad i=1,2,3
$$
for all $x\in \dot{\mathbb{R}}^3_+$. The system  \eqref{seasonal-system} is said to be impermanent if it is not permanent.
\end{definition}
\begin{definition}[\cite{hofbauer1987coexistence,Hutson1982,Gyllenberg2020b}]
    The map $\mathcal{P}$ is said to be permanent if there exists a compact forward invariant set $K\subseteq \dot{\mathbb{R}}^3_+$ such that for every $x\in \dot{\mathbb{R}}^3_+$, there exists $m=m(x)>0$ such that $\mathcal{P}^k(x)\in K$ for all $k\geq m$.
\end{definition}

Let $A$ be the ${3\times 3}$ matrix with entries $a_{ij}$ given in the system \eqref{seasonal-system} and set $r_{i}=b_{i} \varphi \omega-\mu_{i}(1-\varphi) \omega$, $i,j=1,2,3$.

\begin{lemma}[Lemma 2.1 in \cite{Hsu-Zhao}]\label{lemma:exitinction}
If $r_i\leq 0$, then $\Psi_i(t,x)\to 0$ as $t\to +\infty$ for all $x$ on the positive $i$-th axis.
\end{lemma}

\begin{lemma}[Theorem 2.3 in \cite{Niu-Wang-Xie}]\label{carrying-simplex}
The Poincar\'{e} map $\mathcal{P}$ admits a carrying simplex $\Sigma_\mathcal{P}$ if $r_i>0$, $i=1,2,3$.
\end{lemma}

\section{The characterisation of permanence}
Denote  the set of all maps taking $\mathbb{R}_+^3$ into itself by $\mathcal{T}(\mathbb{R}_+^3)$ and the set of all Poincar\'{e} maps associated with the time $\omega$ periodic Lotka-Volterra competition model  \eqref{seasonal-system} with a carrying simplex by $\mathrm{CLVS}(3)$. In symbols:
$$
\mathrm{CLVS}(3):=\left\{\begin{array}{l}
    \mathcal{P}\in \mathcal{T}(\mathbb{R}_+^3): \mathcal{P}=\Phi_{\varphi \omega} \circ L,~\Phi{~ \text{is~the solution~map~of}~}\eqref{LV-system}  \\
    \noalign{\smallskip}
      \text{with~} 0<\varphi<1,\omega,\mu_i,b_i,r_i,a_{ij}>0, i,j=1,2,3
\end{array}
\right\}.
$$

\begin{lemma}[\cite{NWX2023}]\label{lemma-PF}
Let $r_{i}>0$. A point $\theta \in \mathbb{R}_{+}^{3}$ is a fixed point of $\mathcal{P}$ if and only if $\hat{\theta}:=\int_{0}^{\varphi \omega} \Phi_{t}(L \theta) d t$ is a nonnegative solution of the linear algebraic system $(Ax^\tau)_i=r_i$ for all $i$ such that $\theta_i>0$.
\end{lemma}

\begin{proposition}[\cite{NWX2023}]\label{lemma-eigenvalue}
Let $\theta \in \mathrm{Fix}(\mathcal{P},\partial \mathbb{R}^3_+)$ and $\hat{\theta}$ be defined in Lemma \ref{lemma-PF}. Then
$$\lambda_{i}(\theta)=\exp\Bigg\{r_i-\sum_{j=1}^3 a_{i j} \hat{\theta}_{j}\Bigg\}
$$ is an eigenvalue of $D\mathcal{P}(\theta)$ for any $i\notin  \mathcal{K}_{\theta}:=\left\{k: \theta_{k}>0\right\}$.
\end{proposition}

Before we study the permanence of the model  \eqref{seasonal-system}, we first recall a classification provided in \cite{NWX2023} for $\mathrm{CLVS}(3)$. On the boundary of its carrying simplex, the Poincar\'e map $\mathcal{P}$ has three axial fixed points $q_{\{1\}}=(q_1,0,0)$, $q_{\{2\}}=(0,q_2,0)$, $q_{\{3\}}=(0,0,q_3)$ besides the trivial fixed point $0$, and it may has a planar fixed point $v_{\{k\}}$ in the interior of the coordinate plane $\Pi_k$.

Two maps $\mathcal{P}, \hat{\mathcal{P}} \in \mathrm{CLVS}(3)$ with the distinct parameters are said to be {\it equivalent relative to the boundary of the carrying simplex} if there exists a permutation $\sigma$ of $\{1,2,3\}$ such that $\mathcal{P}$ has
a fixed point $q_{\{i\}}$ {\rm(}or $v_{\{k\}}${\rm)} if and only if $\hat{\mathcal{P}}$ has a fixed
point $\hat{q}_{\{\sigma(i)\}}$ {\rm(}or $\hat{v}_{\{\sigma(k)\}}${\rm)}, and
further $q_{\{i\}}$ {\rm(}or $v_{\{k\}}${\rm)} has the same hyperbolicity and local dynamics on the carrying simplex as $\hat{q}_{\{\sigma(i)\}}$ {\rm(}or $\hat{v}_{\{\sigma(k)\}}${\rm)}. $\mathcal{P}\in\mathrm{CLVS}(3)$ is said to be {\it stable relative to the boundary of the carrying simplex} if
all the fixed points on $\partial \Sigma_\mathcal{P}$ are hyperbolic. An
equivalence class is stable if each mapping in
it is stable relative to the boundary of the carrying simplex.

It is proved in \cite{NWX2023} that there are a total of 33 stable equivalence classes for $\mathrm{CLVS}(3)$ via the equivalence relative to the boundary of the carrying simplex which are described in terms of inequalities on parameters. Moreover, classes 1--18 have no positive fixed point so that every orbit converges to some fixed point, while classes 19--33 have at least one positive fixed point. See \cite[Appendix A]{NWX2023} for the parameter
conditions and phase portraits on the carrying simplices for the $33$ classes.

\begin{proposition}\label{prop:per}
    The system \eqref{seasonal-system} is permanent if and only if the Poincar\'e map $\mathcal{P}$ is permanent.
\end{proposition}
\begin{proof}
    If system \eqref{seasonal-system} is permanent,  then there are $\delta,M>0$ such that
    \begin{equation}\label{equ:per}
        \delta\leq \liminf_{t\to +\infty}\Psi_i(t, x)\leq \limsup_{t\to +\infty} \Psi_i(t, x) \leq M,\quad i=1,2,3
    \end{equation}
for all $x\in \dot{\mathbb{R}}^3_+$. It follows that
$$
   \frac{\delta}{2}< \liminf_{k\to +\infty}\mathcal{P}^k(x)\leq \limsup_{k\to +\infty} \mathcal{P}^k(x) < 2M,\quad i=1,2,3
$$for all $x\in \dot{\mathbb{R}}^3_+$.  Then  for all $x\in \dot{\mathbb{R}}^3_+$, there exists $m=m(x)>0$ such that $\mathcal{P}^k(x)\in K_1$ for all $k\geq m$, where $K_1=\{x\in \mathbb{R}^3_+: \frac{\delta}{2}\leq x_i\leq 2M,~i=1,2,3\}$. Let
$K:=\bigcup_{k=0}^\infty P^k(K_1)$. By
\cite[Lemma 2.1]{hofbauer1987coexistence}, $K\subseteq \dot{\mathbb{R}}_+^3$ is a compact forward invariant set of $\mathcal{P}$. Since $K_1\subseteq K$, it follows that  $\mathcal{P}$ is permanent.

On the other hand, if $\mathcal{P}$ is permanent, then there is a compact set $K$ such that  for all $x\in \dot{\mathbb{R}}^3_+$, there exists $m=m(x)>0$ such that $\mathcal{P}^k(x)\in K$ for all $k\geq m$. Let $\Tilde{K}=\Psi([0,\omega]\times K)$. Then for $t\geq (m+1)\omega$, there exists $k\geq m$ and $s\in [0, \omega)$ such that  $t=k\omega+s$ and
$$
\Psi(t,x)=\Psi(k\omega+s,x)=\Psi(k\omega+s,k\omega)\circ\Psi(k\omega,x)=\Psi(s,\mathcal{P}^k x)\in \Tilde{K}.
$$
Since $\Tilde{K}$ is compact, there exist $\delta,M>0$ such that $\delta<x_i<M,~i=1,2,3$ for all $x\in \Tilde{K}$, which implies the model \eqref{seasonal-system} is permanent.
\end{proof}

\begin{lemma}[Corollary 3.5 in \cite{Gyllenberg2020b}]\label{lemma:3-permanence}
Let $n=3$. Suppose a map
$$T(x)=(x_1G_1(x),x_2G_2(x),x_3G_3(x))$$
taking $\mathbb{R}^3_+$ into $\mathbb{R}^3_+$ admits a carrying simplex $\Sigma_T$. Let
$\mathcal{E}(T)=\{x\in \mathbb{R}^3_+: T(x)=x\}$ be the set of fixed points of $T$.
If there are $\nu_1,\nu_2,\nu_3>0$ such that
\begin{equation}\label{Sigma-con-01}
    \eta(\hat{x})=\sum^3_{i=1}\nu_i\ln G_i(\hat{x})>0~(resp. <0),\quad\forall \hat{x}\in \mathcal{E}(T)\cap \partial \Sigma,
\end{equation}
then $T$ is permanent (resp. impermanent).
\end{lemma}

\begin{proposition}\label{coro:3-permanence}
Let $\mathcal{P}\in\mathrm{CLVS}(3)$. If there are $\nu_1,\nu_2,\nu_3>0$ such that
\begin{equation}\label{Sigma-con-02}
    \eta(\theta):=\sum^3_{i=1}\nu_i\ln F_i(\theta)>0~(resp. <0),\quad\forall \theta\in \mathrm{Fix}(\mathcal{P}, \partial \Sigma_\mathcal{P}),
\end{equation}
then the model \eqref{seasonal-system} is permanent~(resp. impermanent).

\end{proposition}
\begin{proof}
The conclusion follows from Proposition \ref{prop:per} and Lemma \ref{lemma:3-permanence}.
\end{proof}

\begin{theorem}\label{permanence-3D}
Assume that $\mathcal{P}\in\mathrm{CLVS}(3)$ is stable relative to the boundary of the carrying simplex. Then
\begin{enumerate}[{\rm (i)}]
\item the model \eqref{seasonal-system} is permanent if $\mathcal{P}$ is in classes $29$, $31$, $33$ and class $27$ with $\vartheta>0$, where
$$
\vartheta:=w_{12} w_{23} w_{31}+w_{21}w_{13}w_{32},
$$
with $w_{ij}=r_j-a_{ji}\frac{r_i}{a_{ii}}$,  $i,j=1,2,3$, $i\neq j$;
\item the model \eqref{seasonal-system} is impermanent if it is in classes $1-26$, $28$, $30$, $32$ and class $27$ with $\vartheta<0$.
\end{enumerate}
\end{theorem}
\begin{proof}
(i) Since the proofs for classes $29$, $31$ and $33$ are completely analogous, we only show the conclusion for the class $31$ here. By \cite[Appendix A]{NWX2023}, for each $\mathcal{P}$ in class 31,   $$\mathrm{Fix}(\mathcal{P},\partial \Sigma_\mathcal{P}) =\{q_{\{1\}},q_{\{2\}},q_{\{3\}},v_{\{1\}},v_{\{3\}}\},$$
where $q_{\{i\}}$ is a axial fixed point on the positive $i$-th axis and $v_{\{k\}}$ is a planar fixed point in the interior of the  coordinate plane $\Pi_k=\{x\in\mathbb{R}^3_+:x_k=0,x_i,x_j>0\}$.

By Proposition \ref{coro:3-permanence}, it suffices to prove that
there are real numbers $\nu_1,\nu_2,\nu_3>0$ such that the following inequalities hold:
\begin{subequations}\label{inequalities}
\begin{align}
  \nu_1\ln F_1(q_{\{1\}})+\nu_2\ln F_2(q_{\{1\}})+\nu_3\ln F_3(q_{\{1\}})>0;\label{ineq:aa} \\[2pt]
  \nu_1\ln F_1(q_{\{2\}})+\nu_2\ln F_2(q_{\{2\}})+\nu_3\ln F_3(q_{\{2\}})>0;\label{ineq:bb}\\[2pt]
  \nu_1\ln F_1(q_{\{3\}})+\nu_2\ln F_2(q_{\{3\}})+\nu_3\ln F_3(q_{\{3\}})>0;\label{ineq:cc}\\[2pt]
  \nu_1\ln F_1(v_{\{1\}})+\nu_2\ln F_2(v_{\{1\}})+\nu_3\ln F_3(v_{\{1\}})>0;\label{ineq:dd}\\[2pt]
  \nu_1\ln F_1(v_{\{3\}})+\nu_2\ln F_2(v_{\{3\}})+\nu_3\ln F_3(v_{\{3\}})>0.\label{ineq:ee}
\end{align}
\end{subequations}
Note that for a fixed point $\theta\in \mathrm{Fix}(\mathcal{P},\mathbb{R}_+^3)$, $F_i(\theta)=1$ for all $i\in \mathcal{K}_\theta$. Therefore, we have
$
F_i(q_{\{i\}})=1,~i=1,2,3
$, and 
$$F_2(v_{\{1\}})=F_3(v_{\{1\}})=F_1(v_{\{3\}})=F_2(v_{\{3\}})=1.$$
Let
$$\hat{v}_{\{k\}}:=\int_{0}^{\varphi \omega} \Phi_{t}(L v_{\{k\}}) d t,~  k=1,3.$$
Then $(A\hat{v}_{\{k\}}^\tau)_j=r_j$ for all $j\neq k$ and $\hat{v}_{\{k\}}=0$ by Lemma \ref{lemma-PF}.
Moreover, it is easy to check that
$$(A\hat{v}_{\{k\}}^\tau)_k<r_k~(>r_k) \Leftrightarrow a_{ki}\beta_{ij}+a_{kj}\beta_{ji}<r_k~(>r_k),$$
where $\beta_{ij}=\frac{a_{jj}r_i-a_{ij}r_j}{a_{ii}a_{jj}-a_{ij}a_{ji}}$,  $i,j,k$ are distinct.
It follows from Proposition  \ref{lemma-eigenvalue} and the condition (ii) of class 31 in \cite[Appendix A]{NWX2023} that
$$ a_{12}\beta_{23}+a_{13}\beta_{32}<r_1 \Leftrightarrow  (A\hat{v}_{\{1\}}^\tau)_1<r_1\Leftrightarrow F_1(v_{\{1\}})=e^{r_1-(A\hat{v}_{\{1\}}^\tau)_1}>1$$
and
$$ a_{31}\beta_{12}+a_{32}\beta_{21}<r_3 \Leftrightarrow  (A\hat{v}_{\{3\}}^\tau)_3<r_3\Leftrightarrow F_3(v_{\{3\}})=e^{r_3-(A\hat{v}_{\{3\}}^\tau)_3}>1.$$
Therefore, \eqref{ineq:dd} and \eqref{ineq:ee} hold for any $\nu_1,\nu_2,\nu_3>0$. On the other hand, by Lemma \ref{lemma-PF} one has
$$\hat{q}_{\{i\}}=\int_{0}^{\varphi \omega} \Phi_{t}(L q_{\{i\}}) d t=\frac{r_i}{a_{ii}},\quad i=1,2,3.$$
Since $\gamma_{12},\gamma_{13},\gamma_{21},\gamma_{23}>0$ 
by condition (i) of class 31 in \cite[Appendix A]{NWX2023}, where $\gamma_{ij}=a_{ii}r_j-a_{ji}r_i$, it follows that
$$F_2(q_{\{1\}}),F_3(q_{\{1\}}),F_1(q_{\{2\}}),F_3(q_{\{2\}})>1$$
by Proposition  \ref{lemma-eigenvalue}.
Thus, \eqref{ineq:aa} and \eqref{ineq:bb} hold for any $\nu_1,\nu_2,\nu_3>0$. The inequality \eqref{ineq:cc} can be written as
\begin{equation}\label{inequalities-2}
  \nu_1\ln F_1(q_{\{3\}})+\nu_2\ln F_2(q_{\{3\}})>0.
\end{equation}
Now fix a $\nu_2>0$. It follows from $\gamma_{32}>0$ that $\ln F_2(q_{\{3\}})>0$, so for sufficiently small $\nu_1>0$ one has
\eqref{inequalities-2} holds.  Let $\nu_3=1$.Then such $\nu_1,\nu_2,\nu_3>0$ ensure that  \eqref{ineq:aa}--\eqref{ineq:ee} hold. This proves that each map $\mathcal{P}$ in class $31$ is permanent. Note that $\partial \Sigma_\mathcal{P}$ is a heteroclinic cycle for class $27$ by \cite[Appendix A]{NWX2023}, so $\mathcal{P}$ is permanent according to  \cite[Corollary 3.7]{Gyllenberg2020b} if $\vartheta>0$.

(ii) For each $\mathcal{P}$ in classes 1--26, 28, 30 and 32, there always exists a fixed point on $\partial \Sigma_\mathcal{P}$ which is an attractor on $\Sigma_\mathcal{P}$ (see \cite[Appendix A]{NWX2023}), so it is impermanent. For $\mathcal{P}$ in class $27$ with $\vartheta<0$, the result follows from \cite[Corollary 3.7]{Gyllenberg2020b}.
\end{proof}

\section{Discussion}\label{sec:cycle}
\par This paper presents permanence and impermanence criteria for
the Lotka-Volterra model of three competing species with seasonal succession \eqref{seasonal-system} via the existence of a carrying simplex, which are finitely computable conditions depending on the boundary fixed points. We provide the exact parameter conditions for permanence and impermanence based on the 33 equivalence classes relative to the boundary dynamics provided in \cite{NWX2023}. Specifically, the model \eqref{seasonal-system} whose Poincar\'e map belongs to the classes 29, 31, 33 and class 27 with a repelling heteroclinic cycle in \cite{NWX2023} is permanent, while it is impermanent if its Poincar\'e map belongs to classes 1--26, 28, 30, 32 and class
27 with an attracting heteroclinic cycle.

However, for model \eqref{seasonal-system}, the permanence generally does not imply the global asymptotic stability of the positive fixed point for the Poincar\'e map $\mathcal{P}$. In fact, there might be multiple positive fixed points in these classes.  For example, it is proved in \cite{nwx2023b} that the class 27 can possess multiple positive fixed points. On the other hand, our numerical simulations show that there might be invariant closed curves, on which every orbit are dense corresponding to the quasiperiodic solutions in \eqref{seasonal-system},  in the classes 29 and 31 and the class 27 with a repelling heteroclinic cycle (see Examples \ref{exm-1}--\ref{exm-3}  below).

\begin{example}[Invariant closed curves in class 27]\label{exm-1}
	Set
	$\omega=10$, $\varphi=0.65$, $\mu_1=0.15$, $\mu_2=0.2$, $\mu_3=0.1$, $b_1=0.3$
	, $b_2=0.3$, $b_3=0.25$, $a_{11}=0.2$, $a_{12}=0.35$, $a_{13}=0.2$, $a_{21}=0.1$, $a_{22}=0.2$, $a_{23}=0.3$,
	$a_{31}=0.8$, $a_{32}=0.1$, $a_{33}=0.3$. It is easy to check that the
	system (\ref{seasonal-system}) with such parameters is in class 27 with $\vartheta>0$, which is permanent. The numerical simulations for the solution of the system and the orbit of the Poincar\'e map $\mathcal{P}$ with the initial value $x_0=(0.3,0.4,0.8)$ are shown in Fig.  \ref{fig-c27-2}, which shows that the system admits
	an attracting invariant closed curve on $\Sigma_\mathcal{P}$.
	
	\begin{figure}[h!]
		\centering
		\begin{tabular}{cc}
			\subfigure[The solution of the system (\ref{seasonal-system})]{
				\label{class27-2-a}
				\begin{minipage}[b]{0.42\textwidth}
					\centering
					\includegraphics[width=\textwidth]{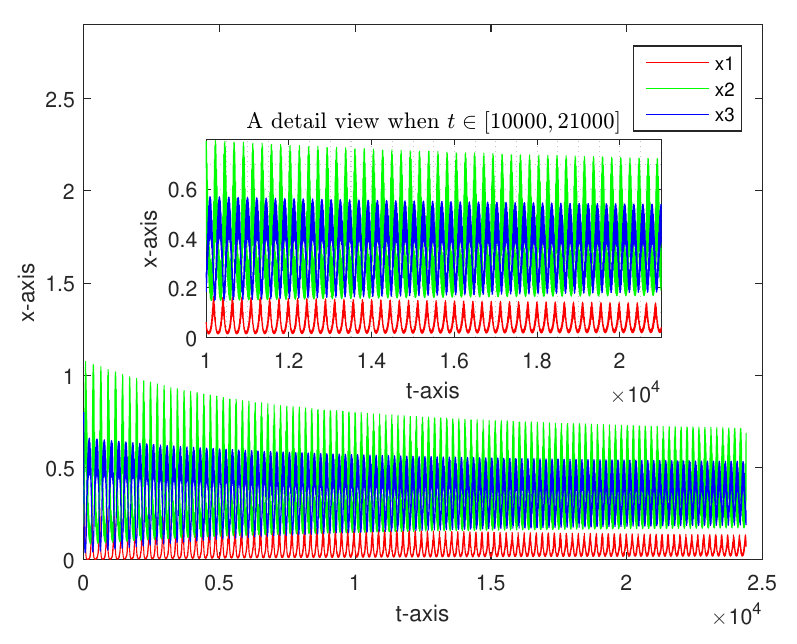}
				\end{minipage}
			} &
			\subfigure[The orbit of the map $\mathcal{P}$]{
				\label{class27-2-b}
				\begin{minipage}[b]{0.42\textwidth}
					\centering
					\includegraphics[width=\textwidth]{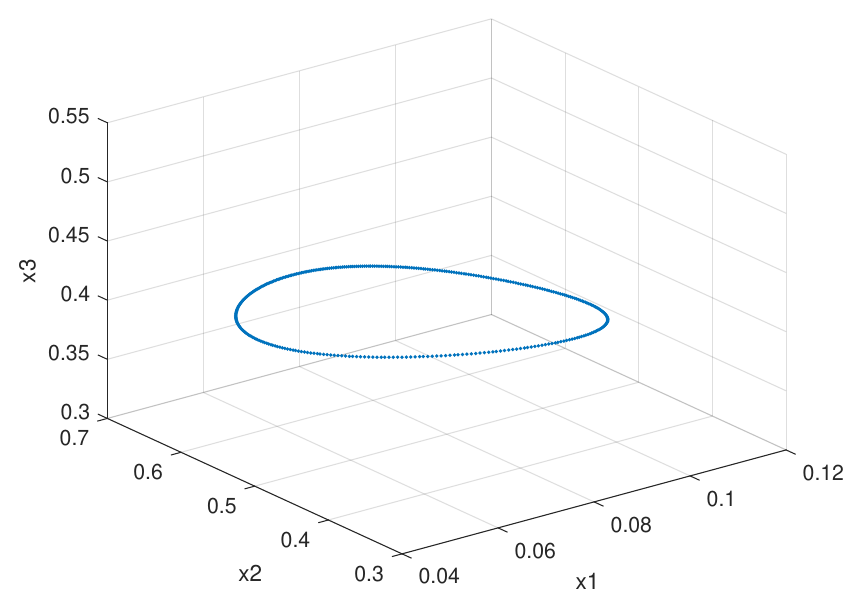}
				\end{minipage}
			} \\
		\end{tabular}
		\caption{The attracting invariant closed curve can occur in class 27.
		}
		\label{fig-c27-2}
	\end{figure}
	
\end{example}

\begin{example}[Invariant closed curves in class 29]\label{exm-2}
	Set
	$\omega=10$, $\varphi=0.5$, $\mu_1=0.1$, $\mu_2=0.11$, $\mu_3=0.15$, $b_1=0.345$
	, $b_2=0.505$, $b_3=0.666$, $a_{11}=0.73$, $a_{12}=0.215$, $a_{13}=0.052$, $a_{21}=1.092$, $a_{22}=0.892$, $a_{23}=0.003$,
	$a_{31}=0.185$, $a_{32}=2.923$, $a_{33}=0.009$. It is easy to check that the
	system \eqref{seasonal-system} with such parameters is in class 29. The numerical simulations for the solution of the system and the orbit of the Poincar\'e map $\mathcal{P}$ with the initial value $x_0=(2,5,2)$ are shown in Fig.  \ref{fig-c29}, which shows that the given system
	admits an attracting invariant closed curve on $\Sigma_\mathcal{P}$.
	
	\begin{figure}[h!]
		\centering
		\begin{tabular}{cc}
			\subfigure[The solution of the system (\ref{seasonal-system})]{
				\label{class29-a}
				\begin{minipage}[b]{0.42\textwidth}
					\centering
					\includegraphics[width=\textwidth]{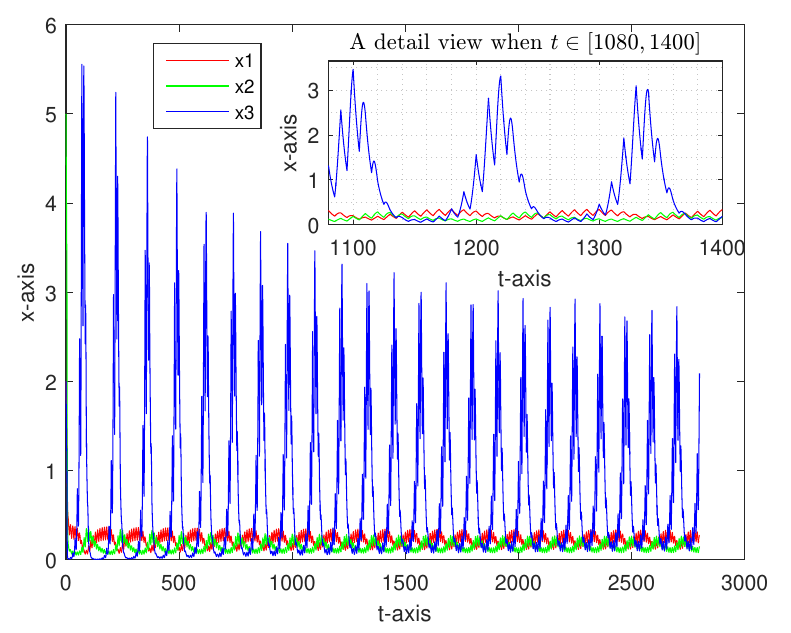}
				\end{minipage}
			} &
			\subfigure[The orbit of the map $\mathcal{P}$]{
				\label{class29-b}
				\begin{minipage}[b]{0.42\textwidth}
					\centering
					\includegraphics[width=\textwidth]{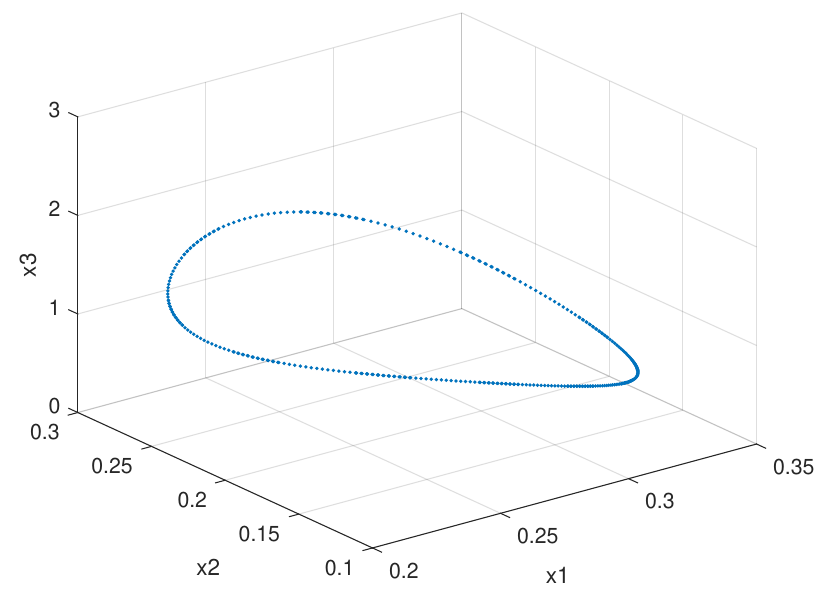}
				\end{minipage}
			} \\
		\end{tabular}
		\caption{The attracting invariant closed curve can occur in class 29.
		}
		\label{fig-c29}
	\end{figure}
	
\end{example}
	
	\begin{figure}[h!]
		\centering
		\begin{tabular}{cc}
			\subfigure[The solution of the system (\ref{seasonal-system})]{
				\label{class31-a}
				\begin{minipage}[b]{0.42\textwidth}
					\centering
					\includegraphics[width=\textwidth]{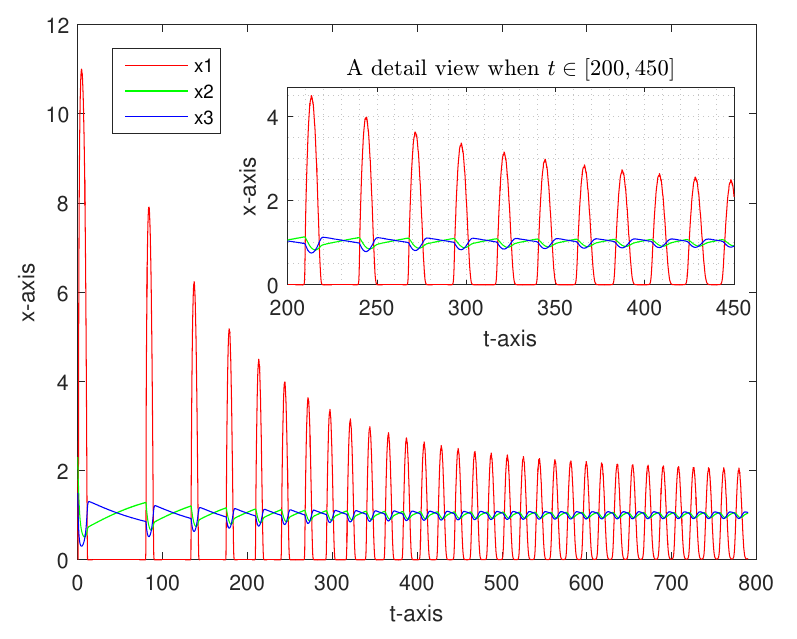}
				\end{minipage}
			} &
			\subfigure[The orbit of the map $\mathcal{P}$]{
				\label{class31-b}
				\begin{minipage}[b]{0.42\textwidth}
					\centering
					\includegraphics[width=\textwidth]{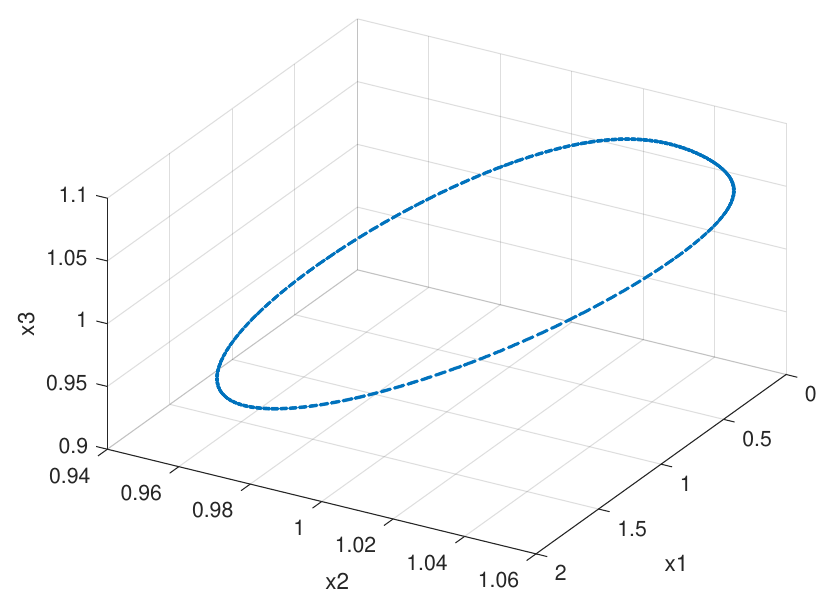}
				\end{minipage}
			} \\
		\end{tabular}
		\caption{The attracting invariant closed curve can occur in class 31.
		}
		\label{fig-c31}
	\end{figure}
\begin{example}[Invariant closed curves in  class 31]\label{exm-3}
	Set
	$\omega=1$, $\varphi=0.97$, $\mu_1=0.23$, $\mu_2=0.27$, $\mu_3=0.18$, $b_1=108$, $b_2=1.2$, $b_3=2.3174$, $a_{11}=6.99$, $a_{12}=1$, $a_{13}=100.2$, $a_{21}=0.074$, $a_{22}=0.521$, $a_{23}=0.602$,
	$a_{31}=0.1174$, $a_{32}=1$, $a_{33}=1.2$. Then the system \eqref{seasonal-system} with such parameters is in class 31. The numerical simulations for the solution of the system
	and the orbit of the Poincar\'e map $\mathcal{P}$ with the initial value $x_0=(1.8,2.3,1.5)$ are shown in Fig.  \ref{fig-c31}, which shows that the system
	admits an attracting invariant closed curve on $\Sigma_\mathcal{P}$.

\end{example}

Interestingly, according to the results in \cite{nwx2023b}, the permanence can indeed imply the uniqueness of the positive fixed point which is globally asymptotically stable for the Poincar\'e map of the system \eqref{seasonal-system} when
all the 3-species are with identical death rates and growth rates, i.e., $\mu_1=\mu_2=\mu_3$ and $b_1=b_2=b_3$. And yet, under which conditions the permanence can guarantee the global stability of the  positive fixed point when $\mu_i$ or $b_i$ are not identical is still unknown. We will leave
it for future research.

 \section*{Acknowledgement}
The authors are very grateful to Prof. Yi Wang for his
valuable and useful discussions and suggestions. Niu’s research was supported by the National Natural Science Foundation of China
(No. 12001096) and the Fundamental Research Funds for the Central Universities (No.
2232023A-02). Xie’s research was supported by the Natural Science Foundation of
Fujian Province (No. 2022J01305).

\bibliographystyle{elsarticle-num-names}
\bibliography{refs}

\end{document}